\documentclass {amsart}

\usepackage{amsfonts,amssymb}
\usepackage{graphicx}
\usepackage{amscd}
\usepackage{eucal}
\input epsf
\newtheorem{thm}{\bf Theorem}[section]
\newtheorem{cor}[thm]{\bf Corollary}
\newtheorem{lem}[thm]{\bf Lemma}
\newtheorem{prop}[thm]{\bf Proposition}

\newtheorem{rem}[thm]{\bf Remark}

\usepackage[dvips]{color}
\definecolor{color1}{RGB}{20,20,20}

\usepackage{amsmath}
\input epsf
\newcommand{\Z}{\mathbb{Z}}

 \newcommand{\W}{\textsf{W}}
\newcommand{\X}{\rm{X}}

\newcommand{\RR}{\textsf{R}}

\newcommand{\V}{\textsf{V}}
\newcommand{\QQ}{\mathbb Q}
\newcommand{\QQQ}{\widetilde {\mathbb Q}}

\def\T{\textsf{T}}
\def\TT{{\mathcal T}}
\def\RRR{{\mathcal R}}
\def\Z{\mathbb{Z}}
\def\A{\textsf{A}}
\def\B{ \textsf{B}}
\def\E{ \textsf{E}}
\def\Q{\textsf{Q}}
\def\SS{\textsf{S}}

\def\V{\textsf{V}}

\def\PSL{{\rm PSL}}
\def\<{\langle }
\def\>{\rangle }
\def\c{{\rm c}}
\def\v{{\bf v}}

\begin{document}


\title{Rational  tangles and  the  modular group}
\author{  Francesca Aicardi}
\maketitle

\begin{abstract} There  is  a natural way  to  associate
with a transformation of an isotopy class  of rational  tangles to  another,    an  element  of  the modular group.   The  correspondence between  the  isotopy classes of rational tangles  and  rational numbers  follows, as well  as the relation with the  braid  group  $B_3$.
\end{abstract}

\section*{Introduction}

The  main result on  rational tangles \cite{Gold} is  the  theorem  stating
that  it  is possible to associate with  every
rational  tangle one and  only one rational  number,  so that  two  rational  tangles are isotopy equivalent iff
they are  represented  by the same  rational  number.

We  obtain  this result  associating with  every   rational tangle an   element  of $\PSL(2,\Z)$.

A   {\it tangle}  with four  ends (here  called simply
tangle) is  an embedding of  two closed  segments  in  a  ball,
such that  their endpoints  are  four  distinct points  of
the   bounding  sphere,  and  the image of  the  interior
of the segments  lie at  the  interior  of the ball.

Let our sphere be centered  at the origin of   the  three space with   Carthesian  coordinates
$X$, $Y$ and $Z$, and     the endpoints of the  strands lie
 in the plane  $Z=0$.

\begin{figure}[h]
\centerline{\epsfbox{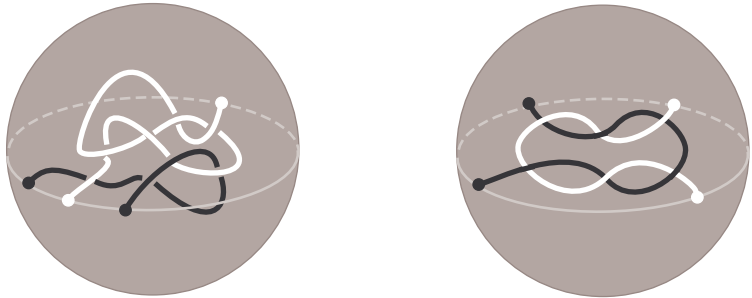}} \caption{Two tangles   }
\end{figure}

We  depict  a  tangle   projected  to  the  $XY$-plane,  which   contains the  endpoints
of the  strands.    Without  lost  of generality, we put  these endpoints
 at  the  intersection of  the  sphere with  the  main  diagonals  of  the plane (Fig. \ref{RK2}).

{\bf Definition.} Two  tangles  are said {\it  isotopic}  if one can  be deformed continuously  to  the other  in the
set of  tangles  with fixed  endpoints.

\section{Rational Tangles}

A tangle  is said  {\it rational} if  can  be  deformed continuously, in the  set  of tangles  with  non  fixed endpoints, to a tangle  consisting  of  two   unlinked  and  unknotted segments.

{\it Example.} The  tangles  in  Figure 1  are  not rational,  since the  white strand
of  the  left  tangle  is  knotted,  and the  two strands  of  the  tangle  at right cannot  be
unlinked  if we allow the  endpoints  of the strands  moving on the bounding  sphere.

For  short we  write   {\it r-tangle}  for rational tangle,  and we denote by $\RRR$ the space of rational tangles.

 \begin{figure}[h]
\centerline{\epsfbox{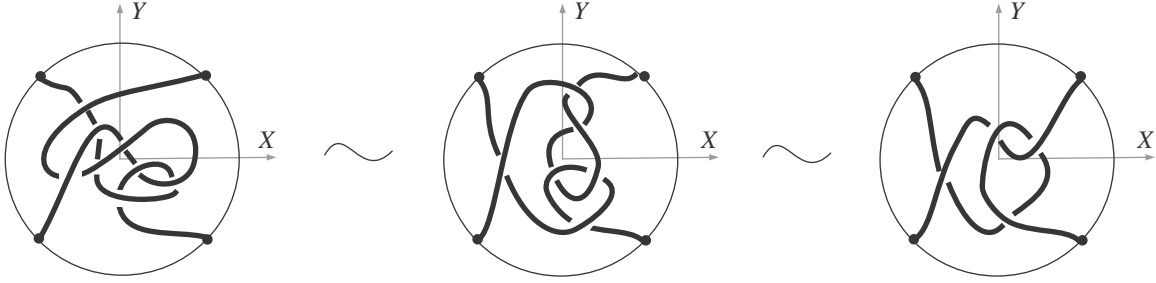}} \caption{Isotopic  r-tangles. }\label{RK2}
\end{figure}

\begin{figure}[h]
\centerline{\epsfbox{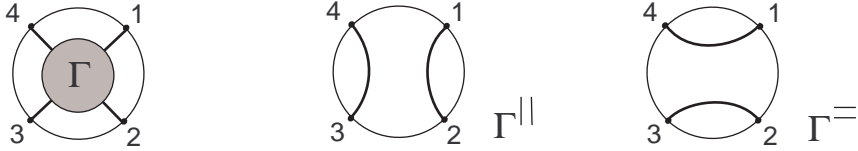}} \caption{Scheme  of a  tangle $\Gamma$  and the basic  tangles $\Gamma^{||}$ and $\Gamma^=$}\label{RK8}
\end{figure}

We denote by $\Gamma^{||}$   and  $\Gamma^{=}$
the   simplest  r-tangles (see  Fig. \ref{RK8})  and   by  1,2,3,4  the  endpoints of the strands,
starting with  point  1 ($X>0, Y>0$) and going  to 4  in  the clockwise  direction.

\begin{figure}[h]
\centerline{\epsfbox{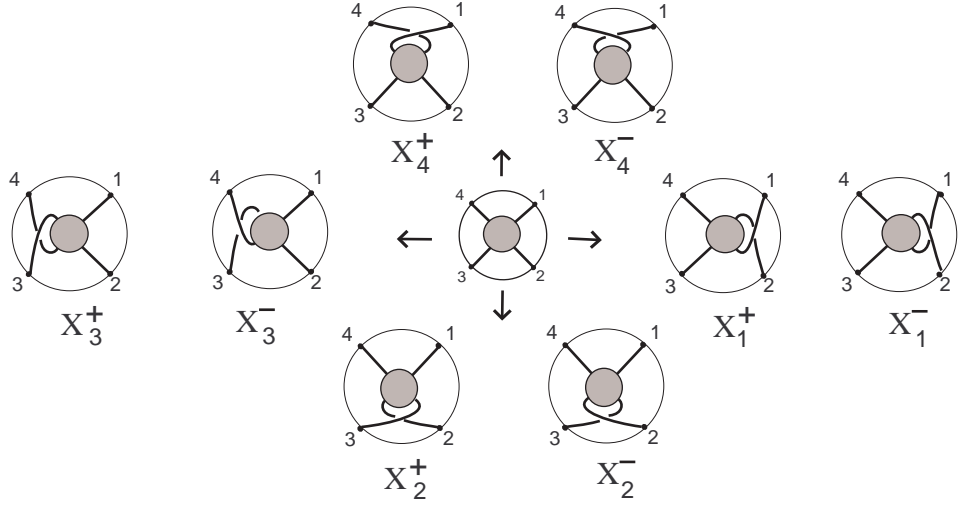}} \caption{The moves  transforming  $\Gamma$  into  $\X_i^\pm \Gamma$}\label{RK10}
\end{figure}

By the constructive definition  of  r-tangles in  \cite{Gold},  an  r-tangle is  obtained
from  $\Gamma^{||}$ or $\Gamma^=$  by a series  of  moves,  consisting in twisting  pairs of  adjacent  endpoints.  We  denote by $\X_i^+$  and   $\X_i^-$  the  {\it  positive} and the  {\it negative}  twists, as shown in Fig.  \ref{RK10}.

Given  an r-tangle  $\Gamma$, the  r-tangle obtained   by applying the  move  $\X_i^\sigma$ ($\sigma=+\  \hbox{ or }\ -$) to $\Gamma$    is  denoted  by $\X_i^\sigma\Gamma$. Thus  any expression
\begin{equation}\label{ex1}  \X_{i_1}^{\sigma_1} \X_{i_2}^{\sigma_2}  \cdots \X_{i_n}^{\sigma_n} \Gamma_0, \end{equation}
where $\Gamma_0=\Gamma^{||}$ or $\Gamma^=$, represents  the r-tangle  obtained from $\Gamma_0$  by  applying
the  moves $\X_i^\pm$ in the order  from $\X_{i_n}^{\sigma_n}$  to $X_{i_1}^{\sigma_1}$.

\begin{figure}[h]
\centerline{\epsfbox{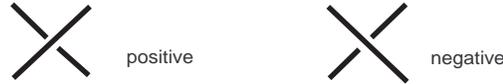}} \caption{Positive  and  negative double points  in  standard representations
of  r-tangles }\label{RK12}
\end{figure}

Any  r-tangle,  written  in  terms  of  the  $\X_i^\pm$,  can be  therefore represented  such  that each   double  point  corresponds to one of the  moves $\X_i^\pm$  and  its sign   is  {\it  positive}  or  {\it negative}
according to the  scheme of  Figure \ref{RK12}. We  call  such   representations    {\it standard} representations.

{\bf  Remark.}  An r-tangle  in a standard  representation, where  the  double  points are  either  all
positive or  all  negative,  is {\it alternating},  i.e., each  strand passes  alternatively  up  and down
in  the  sequence  of double  points encountered  along  it.

Isotopic  r-tangles may  have different standard representations.

{\it Example.}   The representations  of isotopic r-tangles  of  Figure  \ref{RK2},  at  middle  and  at right, are  standard, and they   are respectively:
\[  \X_3^+ \X_2^+ \X_4^- \X_4^- \X_3^+\X_1^+ \Gamma^{=} \ ,  \quad \quad \X_3^+ \X_2^-  \X_3^+\X_1^+ \Gamma^{=}.   \]

\subsection{Symmetries  of the  isotopy classes of  r-tangles}

For  the classification  of the isotopy  classes  of  the  r-tangle  the  following
observation  is essential.

\begin{lem} Every r-tangle $\Gamma$  is isotopic to  the  tangle  $\widehat\Gamma$ obtained from  $\Gamma$
by  a rotation  about  the  axis  $X=0$  and   to  the  tangle $ \Gamma \rangle $ obtained  from  $\Gamma$
by  a rotation  about  the  axis  $Y=0$. \label{pro}
\end{lem}

{\it  Proof.}  The  basic  tangles   $\Gamma^{||}$  and  $\Gamma^{=}$  are invariant under
these   rotations.   The  tangles with  only  one  double point, obtained  as  $\X_1^\pm \Gamma^{=}$, $\X_3^\pm \Gamma^{=}$ or as $\X_2^\pm \Gamma^{||}$, $\X_4^\pm \Gamma^{||}$,  are  also  invariant under such
rotations.  Given  a tangle $\Gamma$, we  write its  expression (\ref{ex1})
in terms  of  $n$  moves.    Therefore we have
\[    \Gamma=\X^{\sigma_1}_{i_1}\Phi,  \]
   where the r-tangle $\Phi$ is expressed  in  terms  of $n-1$  moves.
 Let  us suppose  that the  isotopy class of $\Phi$  is invariant  under  the  considered  rotations, i.e.
 \[  \widehat \Phi  \sim  \Phi  \rangle \sim  \Phi. \]
The  following   figure  shows  that   $\Gamma$  is  isotopy  equivalent,  by consequence,  to  $\Gamma  \rangle$  and  to $\widehat \Gamma$,  when  $\Gamma=\X_1^+\Phi$. Indeed:
  \[ \Gamma  \rangle  =  \X_1^+  \Phi \rangle \sim  \X_1^+\Phi =\Gamma, \]
and
\[    \widehat \Gamma =  \X_3^+ \widehat \Phi  \sim  \X_1^+  \widehat \Phi  \rangle \sim  \X_1^+ \Phi=\Gamma. \]

\centerline{\epsfbox{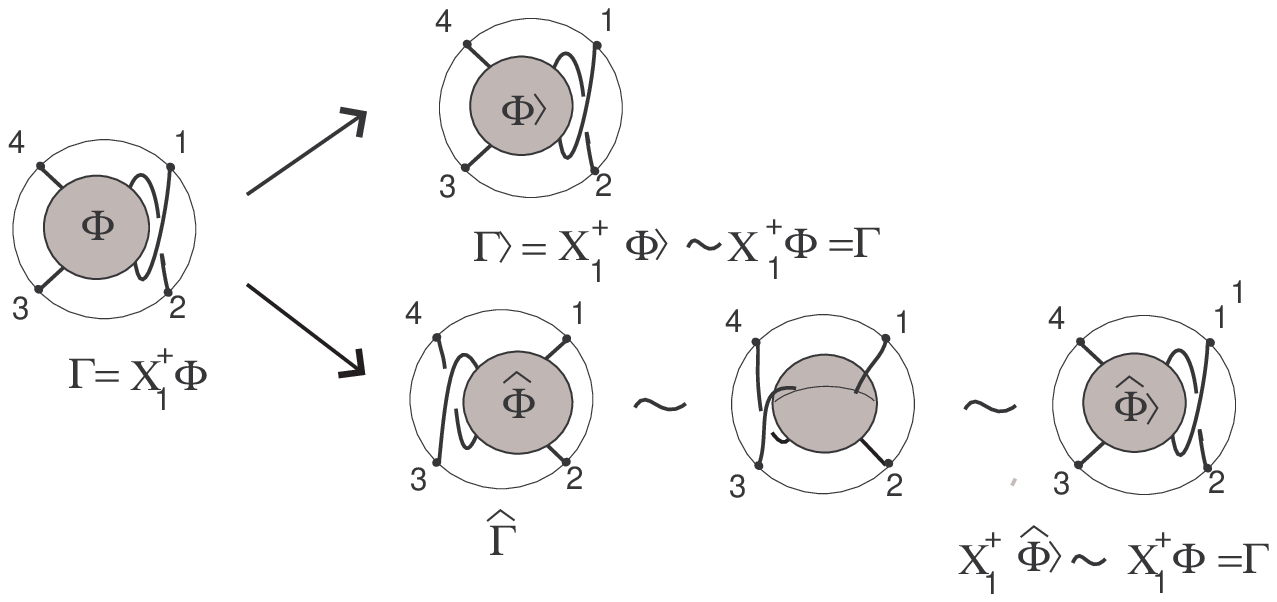}}
Similarly, if  $\Gamma=\X_2^+\Phi$:
  \[   \Gamma\rangle  =  \X_4^+  \Phi\rangle \sim  \X_2^+ \widehat \Phi \rangle  \sim  \X_2^+ \Phi =\Gamma, \]
and
\[    \widehat \Gamma =  \X_2^+ \widehat \Phi  \sim  \X_2^+   \Phi=\Gamma. \]
The  cases with  negative  twists, as  well as the  cases  with $\Gamma=X_3^{\pm}\Phi$  and  $\Gamma=X_4^{\pm}\Phi$   are analogous.   The  lemma  follows.  \hfill $\square$

\section{The  group  of  moves of r-tangles}

From the preceding section  we  deduce evidently    also the following
\begin{cor}\label{XX}   For  every r-tangle  $\Gamma\in \RRR$:
\[  \X_1^\pm  \Gamma  \sim  \X_3^\pm \Gamma, \quad  \quad  \X_2^\pm \Gamma \sim  \X_3^\pm  \Gamma. \]
\end{cor}


Therefore, in order  to classify  the    r-tangles up  to isotopies,
it  is possible  to   consider only  two basic   moves,  that change  the  isotopy  class  of  a r-tangle,
that we  call $\A$ and $\B$,  and precisely:
\[  \A:=  \X_3^+,   \quad  \B:= \X_2^+.  \]
Their  inverses are  denoted,  respectively,  $\A^{-1}$  and  $\B^{-1}$ and  satisfy
$ \A^{-1}=  \X_3^-,   \quad  \B^{-1}= \X_2^-$.

{\it Example.}   The  r-tangles  of  Figure  2  at  middle and at  right,  have   the following  representation:
\[  \A^1 \B^{1}  \B^{-2} \A^2 \Gamma^{=},     \quad \quad \A^1  \B^{-1} \A^2 \Gamma^{=}.  \]

Note that  the moves  $\A$  and  $\B$    keep  invariant  the  tangle endpoint   lying in the  first  quadrant  of the $XY$-plane.

\begin{figure}[h]
\centerline{\epsfbox{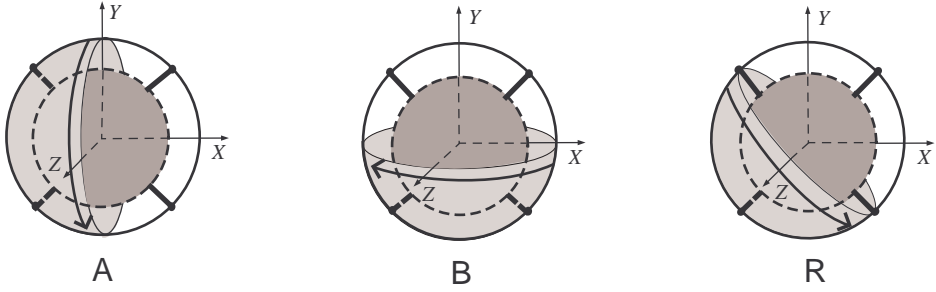}} \caption{The  moves $\A$, $\B$  and $\RR$ obtained  rotating  by $\pi$
one half  of  the  sphere  containing the end-points  of the r-tangle}\label{RK3}
\end{figure}

There  is another  move  keeping  invariant  this  endpoint.  It is the twist  of the endpoints  2 and 4 by a  rotation  by $\pi$    about the  diagonal  $X=Y$, as  shown in  Fig.  \ref{RK3}.  We  denote this move  by $\RR$.

\begin{figure}[h]
\centerline{\epsfbox{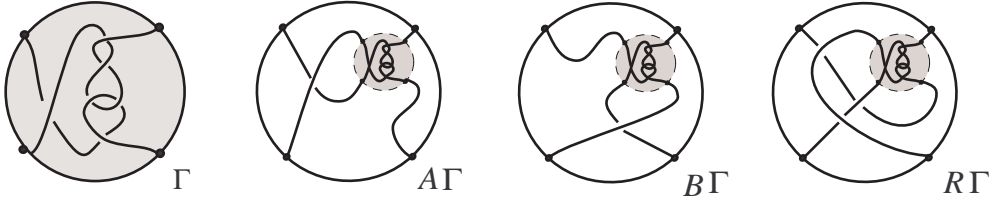}} \caption{Action of the   moves $\A$, $\B$ and $\RR$  on  a  tangle  $\Gamma$   }
\end{figure}

We  now define a  group  of  moves,  generated  by   $\A$, $\B$,  $\A^{-1}$  and $\B^{-1}$,  and we denote it  by  $\TT$.

The  identity   $\E$  is  the   absence  of any  move.

Any  word  $\W=\T_1 \T_2 \cdots \T_n$   where  $\T_i$ are  either   $\A$ or  $\B$, or $\A^{-1}$, or $\B^{-1}$,  represents  an  element of the group, and  will be interpreted as  the  composition  of the moves  $\T_i$, which  are applied   to  any r-tangle  $\Gamma$ starting from  $\T_n$ and ending with  $\T_1$:
\[  \T_1 \T_2 \cdots \T_n \Gamma =  \T_1 ( \T_2 (\cdots (\T_n \Gamma)))).   \]
Two  words  $\W$  and $\Q$  represents  the same  element  iff, for   every $\Gamma\in \RRR$,
\[    \Q \Gamma  \sim  \W  \Gamma. \]
For  instance, the equation  $\A\A^{-1}=\E$  means  that,  for  any  r-tangle $\Gamma$,  $\A \A^{-1}\Gamma \sim  \E \Gamma=\Gamma$.
This equivalence is indeed  obtained  by  the  Reidemeister  move which   eliminates two  double points.

A  sequence  of  $k$ consecutive $\A$ ($\B$)   moves   will be denoted by $\A^k$  ($\B^k$)  and its  inverse by $\A^{-k}$  ($\B^{-k}$).

Any  word $\W=\T_1 \T_2 \cdots \T_n$  will  be thus  written as
\begin{equation}\label{ww} \W=\A^{a_1}\B^{a_2}\A^{a_3}\B^{a_4}\dots  \quad  \hbox{or}  \quad
\W=\B^{a_1}\A^{a_2}\B^{a_3}\A^{a_4}\dots \end{equation}   where $a_i$ are  non zero integers.  The  word $\W$ has  the inverse  \[ \W^{-1}=\dots\B^{-a_4}\A^{-a_3}\B^{-a_2}\A^{-a_1}\quad  \hbox{or}  \quad  \W=\dots\A^{-a_4}\B^{-a_3}\A^{-a_2}\B^{-a_1}\]  satisfying evidently  $\Q^{-1}\Q=\Q\Q^{-1}=\E$.

\begin{thm}  \label{the1}   The
group    $\TT$  is isomorphic  to  $\PSL(2,\Z)$.  \end{thm}

{\it  Proof.} We  define a map   $\mu: \TT\mapsto \PSL(2,\Z)$,     sending  $\A$ and  $\B$    to the following  generators  $A,B$  of  ${\rm PSL}(2,\mathbb{Z})$, which  are defined  up to  multiplications by $- 1$:
\begin{equation}\label{gen}
 A=\left( \begin{array}{cc} 1 & 1 \\ 0 & 1
\end{array}\right), \ \ \ \
 B= \left( \begin{array}{cc} 1& 0 \\ 1 & 1
 \end{array}\right),
 \end{equation}
so that  $\mu(\A^{-1})=A^{-1}$, $\mu(\B^{-1})=B^{-1}$, and    $\mu(\E)=E$, $E$ being the identity $(2\times  2)$-matrix, up  to multiplication  by $- 1$.

The map $\mu$ sends  any sequence  of moves $\A$ and  $\B$ and their inverses  to the  operator of  $\PSL(2,\Z)$
obtained  as the  corresponding product  of $A$, $B$ and  their inverses.

To  prove that $\mu$ is  an  isomorphism,  we have to prove  that
\begin{equation}\label{iso} \mu(\W)=\mu(\W')  \Leftrightarrow \quad  \forall \Gamma \in  \RRR \quad   \W \Gamma \sim \W'\Gamma. \end{equation}
{\it Proof  of} (\ref{iso})$\Rightarrow$.   To every    element of  $\PSL(2,\Z)$  there correspond
different words in  the generators $A$,  $B$, $A^{-1}$   and $B^{-1}$,  since they are not independent.   We  have therefore  to prove that  whenever   $\mu(\W)=\mu(\W')$,
where  $\W$  and  $\W'$  are different words  in  $\TT$, then  $\W\Gamma \sim \W' \Gamma$  for  every r-tangle  $\Gamma$.
This  is   true  if     the following relation  among  the
 considered generators of $\PSL(2,\Z)$  (which  is the  unique non  trivial  relation among them):
\begin{equation}\label{rel} AB^{-1} AB^{-1}A B^{-1}= E \end{equation}
is  the image by  $\mu$  of the  the analog  relation holding in $\TT$:
\begin{equation}\label{rel2} \A\B^{-1} \A\B^{-1}\A \B^{-1}=\E. \end{equation}
Let  $S:=A^{-1} B  A^{-1}
  =\left( \begin{array}{cc} 0 & 1 \\ -1 & 0 \end{array}\right)$, up  to  multiplication  by $-1$.
Observe  that  $S^{-1}= S$  in $\PSL(2,\Z)$. Therefore Eq. (\ref{rel})
is equivalent  to  the  equations
\begin{equation}\label{rel1} B^{-1}A B^{-1}= S, \quad\hbox{and } \quad       S^2=E.  \end{equation}

We  define $\SS:= \A\B^{-1} \A \in  \TT$,  and  to  prove Eq. (\ref{rel2}) we prove:
\begin{enumerate}
\item[(i)] $\B^{-1}\A \B^{-1}= \SS$;
\item[(ii)] $\SS^2=\E$.
\end{enumerate}
Proof of (i).  The next figure shows that for  any  tangle $\Gamma$,
\[ \A\B^{-1} \A\Gamma \sim  \RR\Gamma,  \quad    \B^{-1}\A \B^{-1}\Gamma \sim \RR\Gamma.\]
Hence  item (i) is  satisfied  by  $\SS=\RR$.

  \centerline{\epsfbox{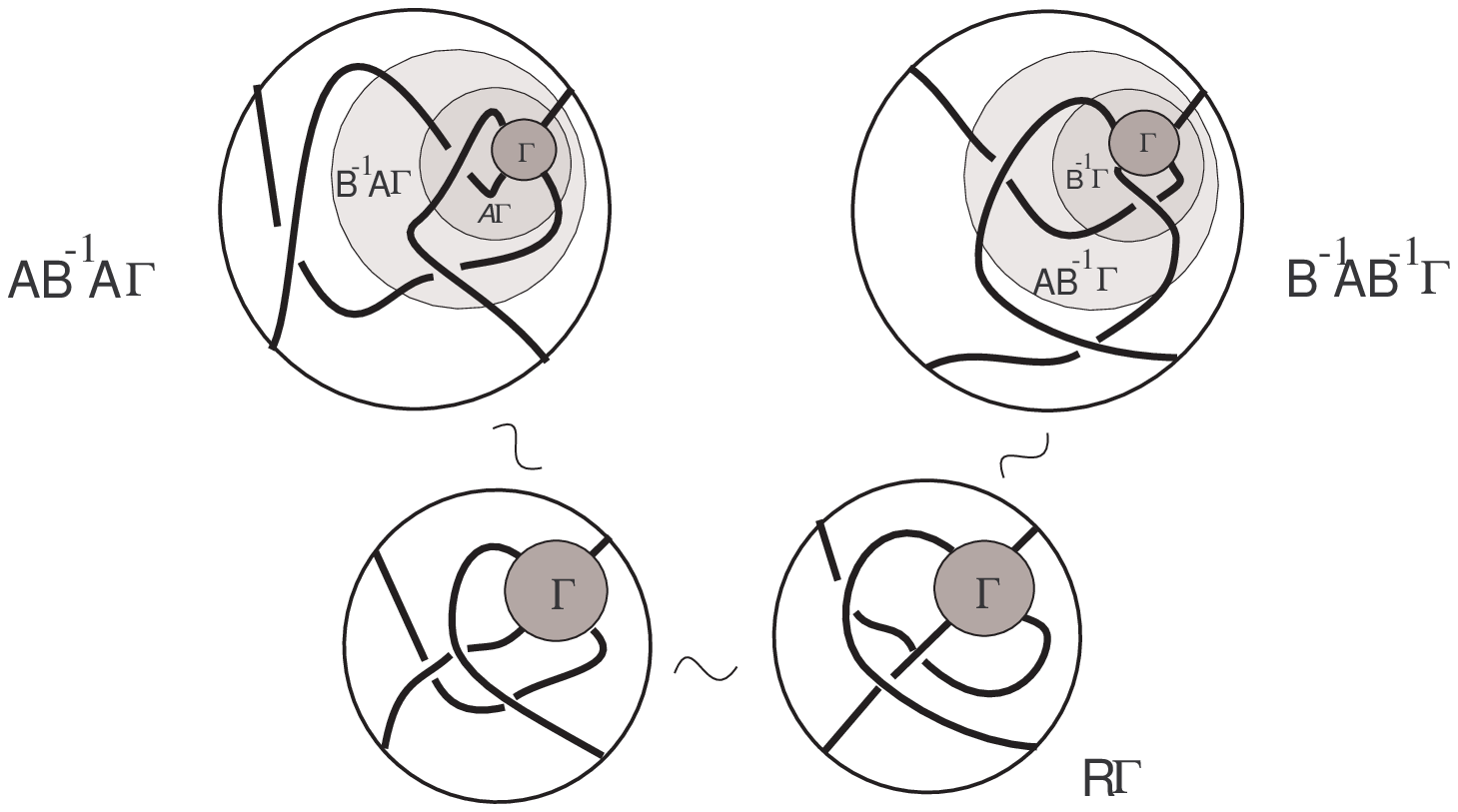}}

Proof of (ii).

\centerline{\epsfbox{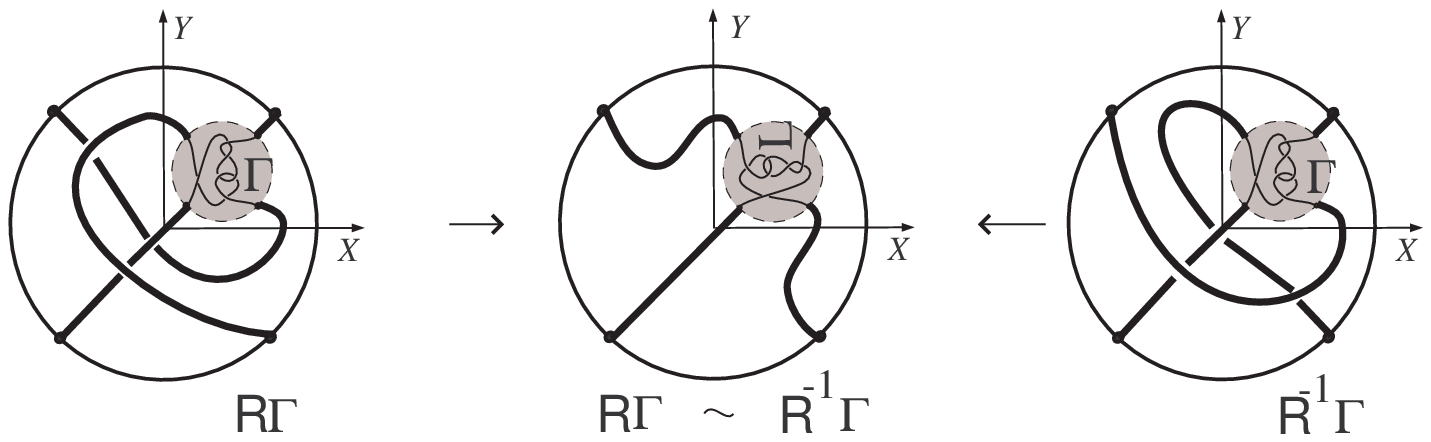}}

Both  $\RR\Gamma$  and $\RR^{-1}\Gamma$  are isotopic  to  the tangle  obtained by
  a rotation of the  whole ball containing $\Gamma$
  about the diagonal $Y=X$.  Hence  $\RR^2  \Gamma\sim \Gamma$
for any  tangle  $\Gamma$, and we conclude $\RR^2=\E$.

{\it  Proof of} (\ref{iso})$\Leftarrow$.
To conclude the proof of the  theorem,  we have  to  exclude that  there  exist
some words $\W$  and  $\W'$    in $\TT$ such that,  for  every r-tangle $\Gamma$,   the isotopy classes  of $\W\Gamma$  and $\W'\Gamma$ coincide,  but  $\mu(\W)\not= \mu(\W')$.   In this case there  should be
  a  move $\Q$, namely  $\Q=\W^{-1}\W'$,  which  is  isotopy  equivalent  to  the
  identity, and   satisfies   $\mu(\Q)\not=E$.
\begin{lem}\label{VT} Every operator  $Q\in \PSL(2,\Z)$  can be  written as  $V T $  where  $T$  is  either  a word
in sole  $A,B$ or  a  word  in sole $A^{-1}, B^{-1}$,  and $V=E$  or $V=S$. \end{lem}

{\it Proof  of the Lemma.} It follows   from the relations (\ref{rel1})   and from  the  derived  relations  holding in  $\PSL(2,\Z)$:
\begin{equation}\label{SASB}    SA=B^{-1}S,  \quad  SB=A^{-1}S,     \quad  SA^{-1}=BS,  \quad  SB^{-1}=AS.  \end{equation}
\hfill  $\square$

We then write the operator $Q:=\mu(\Q)$  in  the  form $Q=VT$.    $T$ is a  word either in  $A,B$
or  in  $A^{-1}, B^{-1}$.  Suppose  that   the word  $T$  ends  with   $A^a$,  for some  non  zero $a\in \Z$. Consider the  r-tangle  $\Q'\Gamma^=$,  where  $\Q'=\V\T$,  $\T$ being  the word  in  the  moves  $\A$
and  $\B$  obtained  translating   $A$  to  $\A$ and $B$   to $\B$ from  the word  $T$,  and  $\V$ is equal  to $\RR$ or to $\E$,  according  to $V$.
The  double  points   of  $\T\Gamma$ are  either all  positive or  all  negative,  and,  since the  move $\RR$  is  equivalent to a rotation  changing the  sign of  all double points, the r-tangle $\Q'\Gamma$ is alternating. It is therefore  evident  that $\Q'\Gamma^= \nsim  \Gamma^=$. But
we have  $\Q'\Gamma^=\sim  \Q\Gamma^=$, since  $\mu(\Q')=\mu(\Q)$.  Therefore $\Q \Gamma^= \nsim  \Gamma^=$.  This contradicts  the  hypothesis  that   the move $\Q$  is  isotopy equivalent to the identity.
If $\T$ ends  with  $\B^a$     for some  non  zero $a\in \Z$,  then  we  consider  the  analog  move $\Q'$  and
  the tangle $\Q'\Gamma^{||}$,  getting  the relations  $\Q \Gamma^{||} \nsim  \Gamma^{||}$ and  $\Q'\Gamma^{||}\sim  \Q\Gamma^{||}$,  from which $\Q \Gamma^{||} \nsim  \Gamma^{||}$, contradicting the hypothesis.
This  concludes  the  proof  of the theorem.  \hfill $\square$

Let  us  denote $\Gamma_0$  the  tangle  $\Gamma^=$.  Observe  that
\begin{equation}\label{AB}
 \Gamma^{||}\sim \B^{-1}\A  \Gamma_0,  \quad \hbox {or} \quad  \Gamma^{||}\sim \A^{-1}\B  \Gamma_0 \end{equation}

We  obtain   the  following corollary:

\begin{cor}\label{space}  The  space  of the isotopy classes  of  rational tangles  coincides with the  orbit of  $\Gamma_0$  under the  action of the  group  $\TT$.  \end{cor}

{\it Proof.}  Every  r-tangle  is isotopy  equivalent  to  an r-tangle reached  from  $\Gamma_0$  or $\Gamma^{||}$
   by  a  series of twists $X_i$.   By  the corollary  \ref{XX},  to  generate  all  the isotopy classes in  $\RRR$
   the moves $\A$,  $\B$ and  their inverses are  sufficient.  The  same $\Gamma^{||}$ can  be  obtained  from  $\Gamma_0$  by  $\A$,  $\B$ and  their inverses, by Eq. (\ref{AB}).  Since every word  in $\A$,  $\B$ and  their inverses  is  an element of  $\TT$,  and  all classes   are  representable starting  from  $\Gamma_0$, the  corollary follows. \hfill $\square$

\section{Rational tangles  and  rational numbers}

We  denote by  $\widetilde \Gamma$ the  isotopy class  of  the  r-tangle $\Gamma$,  and  by  $\widetilde \RRR$
the space of the isotopy classes  of  r-tangles.

In  this  section  we prove  that every  class  of r-tangles different  from  that of   $\Gamma^{||}$ is uniquely  represented by a  rational number,  and,  vice versa,   that  every  rational  number  represents uniquely  a class of rational tangles.

Let
\[ \QQQ:=\{ (p,q) \in  \Z^2 \setminus (0,0)  \ \ \hbox {modulo} \ \  (p,q)\sim (rp,rq), \forall r\in  \Z\setminus 0    \}. \]
Let   $\v_0=(0,1)$  and  $\v_{\infty}=(1,0)$.  Observe  that   $(0,r)$,  for all nonzero $r\in \Z$,  represent  the same element as  $\v_0$  in $\QQQ$,  and $(r,0)$,  for  all  nonzero  $r\in \Z$,  
represent  the  same  element as  $\v_\infty$.

\begin{thm} \label{tm} The  map
$\rho: \widetilde \RRR \rightarrow   \QQQ$   defined  in  this  way:
\[  \rho(\widetilde{\Gamma_0})= \v_0; \]
if  $\Gamma=\Q \Gamma_0$,  then
\[   \rho(\widetilde \Gamma)=\mu(\Q) \v_0,  \]
associates   with   every isotopy  class of  r-tangles  one and only  one  element of  $\QQQ$,
\end{thm}

Before  proving   this  theorem,  we  make some  observations.

\begin{rem}\label{inv} The  r-tangles   $\Gamma_0$  and  $\Gamma^{||}$   possess  the  following
invariances:  for  every  $m\in \Z$
{\rm \[ \B^m \Gamma_0 \sim \Gamma_0 \ ,  \quad  \A^m \Gamma^{||} \sim \Gamma^{||}.
\]  }
No other  element  of  $\TT$  keeps invariant $\Gamma_0$  nor $\Gamma^{||}$.
\end{rem}
Observe  that,  from  Eq. (\ref{AB})  and Remark  \ref{inv} it  follows also  that
\begin{equation}\label{ABR}
  \RR \Gamma_0\sim  \Gamma^{||},  \quad  \RR  \Gamma^{||}\sim  \Gamma_0. \end{equation}

\begin{lem}\label{noinv}    There  is no  element of  $\TT$  keeping  invariant  a  class of  $\widetilde \RRR$
different  from  $\widetilde  {\Gamma_0}$  or  $\widetilde  {\Gamma^{||}}$.
\end{lem}

{\it  Proof.}    Suppose  the  equation $\W\Gamma \sim \Gamma$  be  fulfilled  by  some r-tangle  $\Gamma$ and  some
word  $\W\not= \E$.  Let $\Gamma=\Q\Gamma_0$
for some  $\Q\in  \TT$.      By  hypothesis we  have  $\W\Q\Gamma_0\sim
\Q\Gamma_0$, and, by consequence, $\Q^{-1}\W \Q \Gamma_0 \sim  \Gamma_0$.  By  Remark  \ref{inv}  we obtain that
$\Q^{-1}\W \Q=\B^k$,  for  some integer $k$.  The  last  equation  is   fulfilled either  if both $\W$  and $\Q$  are powers  of $\B$,  and  hence  $\Gamma\sim \Gamma_0$,   and   $\W=\B^k$,  or  if  $\Q=\A^j \RR$,  and  hence,  by  Eq. (\ref{ABR}),  $\Gamma  \sim  \Gamma^{||}$,  and  $\W=\A^{-k}$, being $\RR \A^{-j} \A^{-k} \A^{j} \RR=\B^k$. \hfill $\square$

\begin{lem}\label{TT}   Every  class in  $\widetilde \RRR$   is uniquely represented by  an  element  of    $\TT$,
modulo  multiplication at  right  by {\rm $\hbox{\textsf B}$$^m$, $m\in \Z$}.
\end{lem}

{\it  Proof.}  Every  class  of  r-tangles  can be  represented as $\Q\Gamma_0$,  for  some  $\Q\in \TT$,
by  Corollary   \ref{space}.  With the r-tangle   $\Gamma=\Q \Gamma_0$  we associate the element  $\Q \in \TT$. By  Remark
\ref{inv},  $\widetilde \Gamma=\widetilde{\Q \B^m \Gamma_0}$,  therefore  if we associate $\Q$  to $\Gamma$,
we  may as well  associate  $\Q\B^m$, for all  $m\in \Z$. On  the  other  hand,  if  $\Q'\not=\Q \B^m$ for some  $m$, then $\Q'\Gamma_0\nsim \Q\Gamma_0$,  by the  same  remark.   \hfill  $\square$

{\it  Proof of Theorem} \ref{tm}.  We  observe firstly that
the map $\rho$  respects the  invariances  of Remark \ref{inv}.

Indeed,  $\rho(\widetilde {\B^m \Gamma_0})=B^m\v_0=(^1_m  \ ^0_1) (^0_1)= (^0_m )=m\v_0\sim \v_0=\rho(\widetilde{\Gamma_0})$.
Moreover,  using  Eq. (\ref{AB}),
\begin{equation}\label{roinf} \rho( \widetilde {\Gamma^{||}})= B^{-1}A \v_0=(^{\ 1}_{-1} \ ^1_0)(^0_1 )=(^1_0)=\v_{\infty}.  \end{equation}

Therefore  $\rho(\widetilde {\A^m \Gamma^{||}})=A^m B^{-1}A \v_0=A^m\v_\infty=(^1_0 \ ^m_1) (^1_0)= (^m_0 )=m\v_\infty\sim \v_\infty$,
i.e., $\rho( \widetilde {\A^m \Gamma^{||}})=\rho( \widetilde {\Gamma^{||}})$.

We have  to  prove that  for  every pair  of  classes  $\widetilde \Gamma$  and  $\widetilde{\Gamma'}$  in $\widetilde \RRR$
\begin{equation}\label{teo} \widetilde \Gamma = \widetilde {\Gamma'} \quad  \Leftrightarrow \quad  \rho(\widetilde \Gamma)=\rho(\widetilde {\Gamma'}). \end{equation}

{Proof of} (\ref{teo})$\Rightarrow$. By  Lemma  \ref{TT}, if $\widetilde \Gamma = \widetilde {\Gamma'}$  then
$\Gamma=\Q\Gamma_0$  and $\Gamma'=\Q'\Gamma_0$  with  $\Q'=\Q\B^m$.
If      $\mu(\Q)=\left(\begin{matrix}a&b\\c & d \end{matrix}\right)$, we  obtain   $\mu(\Q')=\left(\begin{matrix}a&b\\c & d \end{matrix}\right) \left(\begin{matrix}1&0\\m & 1 \end{matrix}\right)= \left(\begin{matrix}a+mb&b\\c+md & d \end{matrix}\right)$.
Therefore  $\rho(\widetilde \Gamma)= \mu(\Q)\v_0= (b,d)$  and  $\rho(\widetilde \Gamma')= \mu(\Q')\v_0= (b,d)$.

{Proof of} (\ref{teo})$\Leftarrow$.  Let  $\Gamma\sim \Q\Gamma_0$,  $\Gamma'\sim\Q'\Gamma_0$,  $\mu(\Q)=\left(\begin{matrix}a&b\\c & d \end{matrix}\right)$,   $\mu(\Q') =\left(\begin{matrix}a'&b'\\c' & d' \end{matrix}\right)$.  We suppose  that   $\rho(\widetilde\Gamma)=\rho(\widetilde {\Gamma'})$, i.e.,
 $(b,d)\sim (b', d')\in  \QQQ$.    This equality implies  $b=b'$  and  $d=d'$, because $b,d$ as  well  as   $b',d'$  are relatively prime,  constituting  a column of a discriminant one matrix.   The matrices  $Q:=\mu(\Q)$  and $Q'=\mu(\Q')$ may  differ  for the first column,
and  they have the same unit discriminant.   We obtain therefore  $a'=a+bk$  and   $c'=c+dk$ for  some  integer $k$,
and this  implies  $Q'=Q B^k$.
We  thus  obtain $\Gamma'\sim \Q' \Gamma_0 \sim \Q \B^k\Gamma_0 \sim \Q\Gamma_0$,   and  hence $\Gamma'\sim \Gamma$.   \hfill $\square$

Consider the  map  $t$  sending every  element of $\QQQ$  with $q\not=0$  to   a  rational  number:  $t:(p,q)\mapsto  p/q$.  Every  rational  number is  the image  by  $t$  of one  and   only  one  element  of  $\QQQ$. We  will write $p/q=t(\v)$,  where  $\v=(p,q)\in  \QQQ$.   Evidently $m\v=(mp,mq)\sim \v$  in  $\QQQ$.

\begin{cor} The map $\rho \circ t: \widetilde \RRR \rightarrow  \mathbb Q$ associates with   every isotopy  class of  r-tangles different  from  $\widetilde{\Gamma^{||}}$   one and only  one  rational number.
\end{cor}

{\it  Proof.} By Theorem \ref{tm},  the map  $\rho$    associates  with  every isotopy  class of  r-tangles  an element of  $\QQQ$. The  map $t$  associates with every element  of $\QQQ$, $q\not=0$,  a rational  number.
  The element  $(r,0)\in \QQQ$  is  the  image by  $\rho$  of  the class of  $\Gamma^{||}$,  by  Eq. (\ref{roinf}). The corollary  follows.

\begin{cor} The map $(\rho \circ t)^{-1}:  \mathbb Q \rightarrow  \widetilde \RRR  $ associates with   every rational  number one and only one isotopy  class of  r-tangles.
\end{cor}

{\it  Proof.}  With  the  rational number  $p/q$,  $t^{-1}$ associates  $\v=(p,q)\in \QQQ$.  Since  $q\not=0$,
and  $p$  and $q$ are  coprime,  the pair $(p,q)$  defines a matrix  $Q\in  \PSL(2,Z)$ such that  $Q\v_0=\v$,
up to  a  right factor  equal  to  $B^m$,  for  some  $m\in \Z$  as we  have  seen  in   the  proof  of Eq. (\ref{teo}).
Therefore also  the  element $\Q\in  \TT$  such that  $\mu(\Q)=Q$ is  defined     up  to a  right factor  equal  to $\B^m$, for  some  $m\in \Z$.
But all elements  $\Q\B^m$  define the  same class  $\widetilde {\Q \Gamma_0}$,  therefore  the  image of $\v$  by
$\rho^{-1}$   is  well  defined  and unique  by Remark \ref{inv}.  \hfill $\square$

We  obtain  as  corollary that  every r-tangle is alternating.

\begin{cor} A  rational tangle  different  from  $\Gamma^{||}$ is equivalent to  a  tangle  obtained by $\Gamma_0$ either by sole positive twists
or  sole  negative twists.   \end{cor}

{\it  Proof.} We  write  $\Gamma\sim \Q \Gamma_0$.  Hence  we  consider  $Q=\mu(\Q)$.  $Q$  is  different  from
$A^m S$  since  $\Gamma\nsim \Gamma^{||}$,  according to   Eq. (\ref{AB}), Remark  \ref{inv} and Eq. (\ref{ABR}).
By  Lemma  \ref{VT} and  Theorem \ref{the1},  $\Q$  can  be written  as  $\V\T$, where  $\T$  is a  word
either in  $\A$ and $\B$  or  in  $\A^{-1}$  and $\B^{-1}$,  and  $\V=\RR$  or  $\V=\E$.  If  $\V=\E$,  the
proof is  finished,  since  the  tangle  $\T \Gamma$  is  alternating.   If  $\V=\RR$,  then we apply  to every generator  in the word  $\T$  the  equations  obtained  by  the  relations  in $\TT$   corresponding  to  the relations  (\ref{SASB}) in $\PSL(2,\Z)$.
Therefore  we obtain  $\RR \T=\bar \T  \RR$,  where  $\bar  \T$ is obtained by $\T$  exchanging  $\B$  with  $\A^{-1}$ and  $\A$  with $\B^{-1}$.   Hence  $\Gamma \sim  \bar \T  \RR  \Gamma_0\sim  \bar \T  \Gamma^{||}$.  Since  $\Gamma^{||}\sim \A^m\Gamma^{||}$,  we can  write $\Gamma \sim  \bar \T'  \Gamma^{||}$, where  $\T'$  is  a word
whose  last  element  is  $\B^k$,  with    $k\not=0$.   If  $k\le 1$,  and hence $\T'$  is  a  word  in
$\A^{-1}$  and  $\B^{-1}$,  we  use the  relation,  from  Eq. (\ref{AB}),   $\B^{-1}\Gamma^{||}\sim \A^{-1}\Gamma_0$,
to obtain  the relation $\T'\Gamma^{||}=\T'' \Gamma_0$, where  $\T''$  is  obtained  substituting the last  $\B^{-1}$ with  $\A^{-1}$.
  In  this way $\Gamma \sim \T'' \Gamma_0$  is obtained  from  $\Gamma_0$ by sole negative  twists.    If  $k\ge 1$,  and hence $\T'$  is  a  word  in
$\A$  and  $\B$,  we  use the  relation,  from  Eq. (\ref{AB}),   $\B\Gamma^{||}\sim \A\Gamma_0$,
to obtain  the  relation  $\T'\Gamma^{||} \sim \T'' \Gamma_0$, where  $\T''$  is  obtained  substituting the last  $\B$ with  $\A$.
Therefore  $\Gamma \sim \T'' \Gamma_0$  is obtained  from  $\Gamma_0$ by sole positive  twists.\hfill  $\square$

\subsection{Continued fraction procedure} We  prove  now  that  the  continued  fraction procedure  allows to  associate  with  every rational  number  one  and only one
 alternating  rational tangle.

 Let  $p,q$ and  $a_i$  ($i=1...n$)  represent  natural numbers.
A positive rational number  has  a unique representation  by a  continued fraction  with  positive  elements $a_i$:
\[   \frac{p}{q}=  a_1+\frac{1}{a_2+\frac{1}{a_3+\frac{1}{ a_4+\frac{1}{...+\frac{1}{a_n}} }   }  }.  \]
As  a consequence, a  negative  rational number  has  a  unique representation  by a  continued fraction  with negative elements $-a_i$:
\[  - \frac{p}{q}=  -a_1-\frac{1}{a_2+\frac{1}{a_3+\frac{1}{ a_4+\frac{1}{...+\frac{1}{a_n}} }   }  }= (-a_1)+\frac{1}{(-a_2)+\frac{1}{(-a_3)+\frac{1}{ (-a_4)+\frac{1}{...+\frac{1}{(-a_n)}} }   }  }.  \]

A continued  fraction  with positive elements  is denoted by  \[  p/q=[a_1;a_2,\dots, a_n].\]
Note  that  $a_1$  may  be  zero, but   all $a_i$, for $i>1$,  are  positive,  and   $a_n$  is   bigger or equal  to  2.
Therefore,   we  can  always  write  an  equivalent  expression  for $p/q$:
\[   \frac{p}{q} =[a_1,a_2,a_3,\dots,  a_n-1,1].\]

{\it  Example}: $7/3=[2,3]=[2,2,1]$,  since      \[  2+\frac{1}{3}= 2+\frac{1}{2+\frac{1}{1}}.   \]

If the fraction is  is negative, we will  write  evidently  $-\frac{p}{q} =[-a_1,-a_2,-a_3,\dots,  -a_n+1,-1]$.

In  this way  we  may always  suppose  $n$ be  odd, because if  $n$  is  even,  then  we  decrease $a_{n}$ by one
and we add  the  $(n+1)$-th  element  equal  to  1 (or $-1$).  We call such a  continued  fraction {\it odd} continued  fraction.

\begin{prop}  Let   the  odd continued  fraction of $p/q$  be
\begin{equation}\label{cfrac} [a_1; a_2,a_3,\dots, a_n], \end{equation}  where the  $a_i$
 are either all  positive or all negative  but  $a_1$ that may  be  zero.
Then  {\rm $p/q=t(\rho(\widetilde{\Q\Gamma_0}))$}, where  $\Q\Gamma_0$  is  the alternating  r-tangle
obtained  applying to  $\Gamma_0$  the sequence of moves {\rm
\[    \A^{a_1}\B^{a_2}\A^{a_3} \cdots  \B^{a_{n-1}}\A^{a_n}. \]}
\end{prop}

{\it  Proof.}  
   We define  the  map $\tau_Q: \QQ \rightarrow \QQ$ in  the following way: for every  $z=\frac{p}{q}\in  \QQ$
\[  \tau_Q (z)= t(Q(^p_q)). \]
Evidently, $t(Q(^p_q))=t(Q (^z_1))$.  For  every $\alpha\in  \Z$  we  have:
\[ \tau_{A^{\alpha}}(z)=z+\alpha, \quad  \quad    \tau_{B^\alpha}(z)=\frac {1}{\alpha+\frac{1}{z}}.   \]
If $\alpha=0$,  we obtain  $\tau_E(z)=z$.
Moreover,  if  $Q''=QQ'$, then  $\tau_{Q''}(z)=\tau_Q(\tau_{Q'}(z))$.
Therefore,  if  $Q=A^{a_1}B^{a_2}A^{a_3} \cdots  B^{a_{n-1}}A^{a_n}$  then
\[  \tau_Q(0)=a_1+\frac{1}{a_2+\frac{1}{a_3+\frac{1}{ a_4+\frac{1}{...+\frac{1}{a_n}} }   }  }.   \]
We  thus  associate  with   the odd  continued fraction  (\ref{cfrac}) of $p/q$
the  sequence of  moves  \[ \Q= \A^{a_1}\B^{a_2}\A^{a_3} \cdots  \B^{a_{n-1}}\A^{a_n} \]  and  the  alternating
tangle  $\Q \Gamma_0$.  By  construction,     $p/q=t(\rho(\widetilde{\Q\Gamma_0}))$.  \hfill $\square$

\section{Rational  tangles  and  the braid  group $B_3$}
It  is  well  known  that the modular group is  isomorphic  to  $\widetilde B_3:=B_3/\langle \omega \rangle $ where $B_3$  is the group of braids  with  3  strands,  generated  by  $\sigma_1$  and  $\sigma_2$:

\centerline{\epsfbox{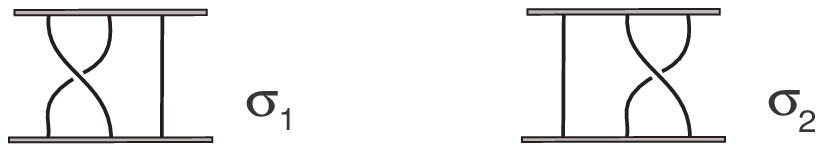}}
and  $\langle \omega \rangle$  is the    relation
\begin{equation} \label{relb} \sigma_1\sigma_2 \sigma_1  \sigma_2 \sigma_1 \sigma_2 = 1.\end{equation}

We describe here  explicitly  the isomorphism  between $\widetilde B_3$ and $\TT$,  giving its
topological meaning.

\begin{thm}\label{corr} The  following     construction  associates with  every  element of $\TT$
an element  of  $\widetilde B_3$  and  vice-versa.
\end{thm}

\centerline{\epsfbox{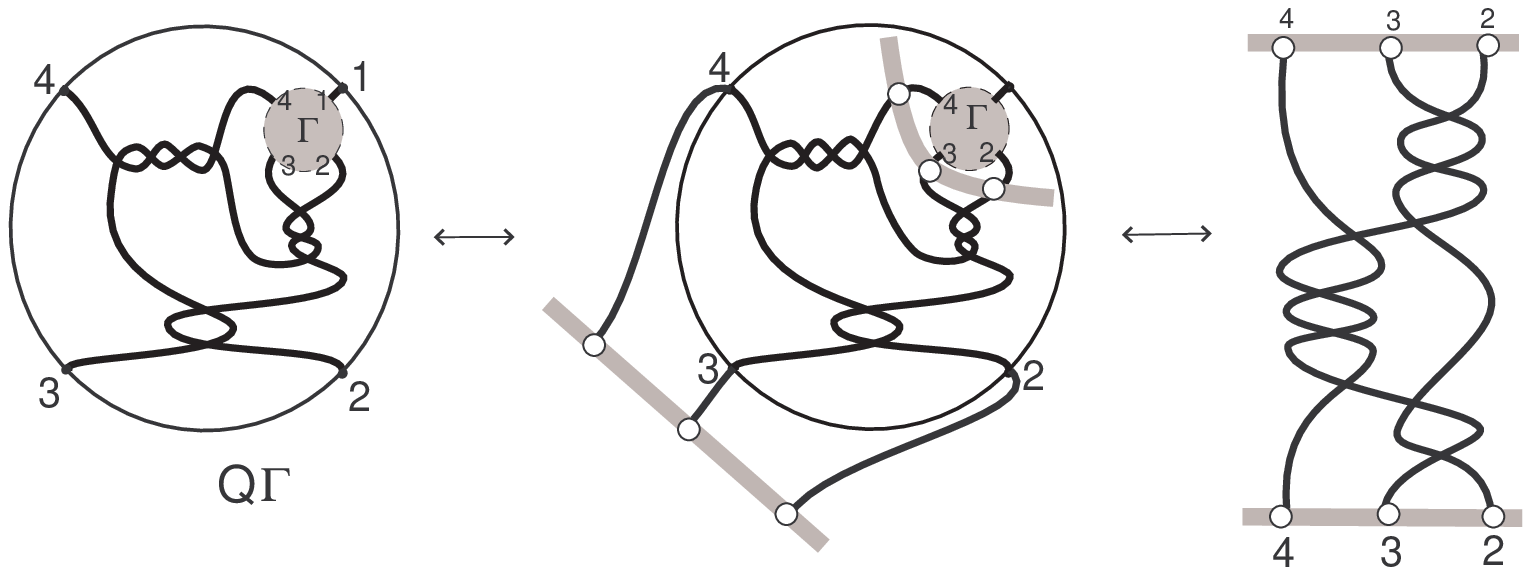}}

{\it  Proof.} The  construction  can  be described as  follows: the  endpoints 4, 3  and  2  of the  r-tangle  $\Gamma$
become, respectively, the  upper endpoints    of  the strands of the  braid,  and the  endpoints 4, 3  and  2  of the  r-tangle  $\Q\Gamma$
become, respectively, the  lower  endpoints  of  the strands of the  braid.

The topological resolution of each double  point remains  unchanged in the transformation,  so that we obtain:
\begin{equation} \label{ABsigma}\A  \longleftrightarrow \sigma_1   \quad \quad  \quad  \B \longleftrightarrow  \sigma_2^{-1}.  \end{equation}

\centerline{\epsfbox{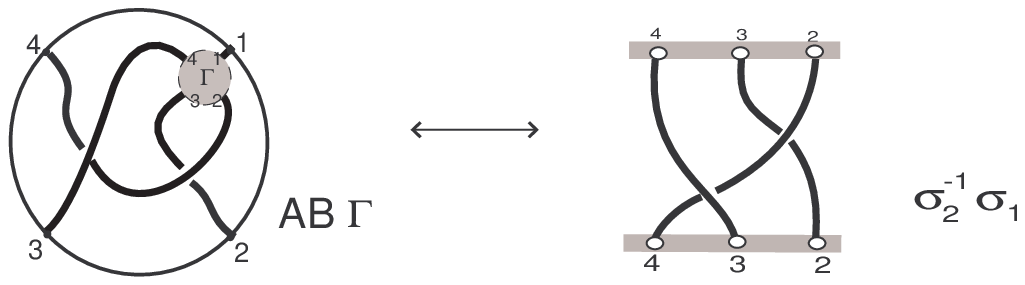}}

To any word of  $\TT$ in $\A$ and $\B$ and  their inverses   there corresponds  therefore  the word in $B_3$
where the  generators are  translated  according to  (\ref{ABsigma})  and  are put in  the inverse order,
if the  order  from left to right of the generators  in  a word  of  $B_3$  is (as usually)  interpreted  as their
order  from  top  to  bottom  in  the  braid.

The relation  $\langle \omega \rangle$  is therefore the  translation  of  relation: $\B^{-1} \A\B^{-1}\A \B^{-1}\A=\E$, which is  equivalent  to (\ref{rel2}).    \hfill  $\square$

{\bf Remark.} Remember that  in  $\RRR$ the  topological meaning  of  relation  (\ref{rel2})  is that  $\RR^2=\E$,  i.e.,  rotating  by $2\pi$  an r-tangle about  the diagonal   containing  the  endpoints 1 and 3, its isotopy class  does  not  change.

For the  braids the  element   $\sigma_1\sigma_2 \sigma_1  \sigma_2 \sigma_1 \sigma_2$  is  as well equivalent  to a  rotation  by $2 \pi$  about the  central vertical axis of  the  segment containing the  lower endpoints (or of that containing  the  upper  endpoints) of the  strands. This rotation  changes  the  isotopy class  of the braid. Therefore $\widetilde B^3$ can  be  interpreted  as  the group  of  braids  such that   the segments  containing  the endpoints   are allowed  to  rotate independently  by $2\pi$.

\centerline{\epsfbox{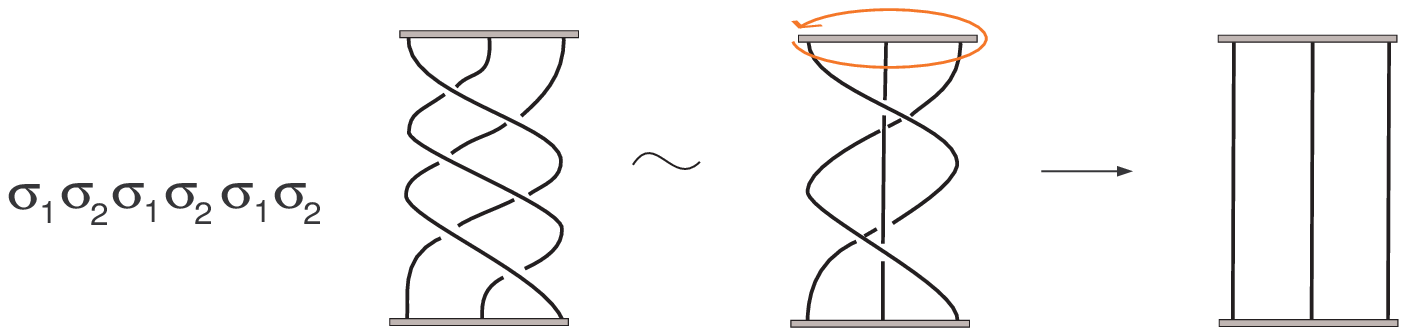}}

{\bf  Remark.} The  construction  of Theorem \ref{corr} is  similar to  that  showed in  \cite{Sof},  but we  haven't any
ambiguity  problem, since  we deal  with  $\TT$,  and  not  with rational tangles.

\end{document}